\documentclass{wsbart/wsbart}

\usepackage{lmodern}
\usepackage{amsmath,amssymb,amsthm,mathtools}
\usepackage{microtype,dsfont}
\usepackage{xurl}
\usepackage[maxnames=99,style=trad-plain,backend=biber,isbn=false,
uniquename=init]{biblatex}
\DeclareRedundantLanguages{English}{english}
\AtBeginBibliography{\setlength{\emergencystretch}{2em}}

\let\authblkAuthor\author
\let\authblkMaketitle\maketitle
\usepackage{titling}
\let\author\authblkAuthor
\let\maketitle\authblkMaketitle
\predate{}
\postdate{}

\usepackage{changes}

\renewcommand{\mathbb}{\mathds}

\renewcommand{\le}{\leqslant}
\renewcommand{\ge}{\geqslant}

\newtheorem{theorem}{Theorem}
\newtheorem{proposition}[theorem]{Proposition}
\newtheorem{lemma}[theorem]{Lemma}
\newtheorem{corollary}[theorem]{Corollary}
\theoremstyle{remark}
\newtheorem{remark}[theorem]{Remark}

\newcommand{\R}{\mathbb R}
\newcommand{\E}{\mathbb E}
\newcommand{\Prob}{\mathbb P}
\newcommand{\dd}{\,\mathrm d}
\newcommand{\one}{\mathds 1}
\newcommand{\BV}{\mathrm{BV}}
\DeclareMathOperator{\Cov}{Cov}
\DeclareMathOperator{\Var}{Var}
\DeclareMathOperator{\tr}{tr}
\DeclareMathOperator{\diver}{div}
\DeclareMathOperator{\Per}{Per}
\DeclareMathOperator{\supp}{supp}
\DeclarePairedDelimiter{\norm}{\lVert}{\rVert}

\title{One-third Hölder stability of optimal transport maps
from log-concave sources}
\author[1]{Songbo Wang}
\author[2]{Xiaozhen Wang}
\affil[1]{CEREMADE, Université Paris-Dauphine, Université PSL,
Paris, France}
\affil[2]{CMAP, École polytechnique, Institut polytechnique de Paris,
Palaiseau, France}
\date{}
\datefootnote{\today}

\begin{document}

\maketitle

\begin{abstract}
We prove quantitative stability of quadratic optimal transport maps
from arbitrary log-concave sources to bounded targets. If the source
$\rho$ is isotropic on $\R^d$ and the targets $\mu$, $\nu$ are supported
in a common ball of radius $R$, their Brenier maps satisfy
\[
 \norm{T_\mu-T_\nu}_{L^2(\rho)}
 \le17d^{2/3}R^{2/3}W_2(\mu,\nu)^{1/3}.
\]
The one-third exponent is optimal for $d\geqslant2$.
The proof combines a linear stability estimate for Brenier potentials
in a weighted norm with a matching weighted reverse Poincaré inequality.
The weight is chosen in terms of the negative log-density of the source
to account for its tails. We also show that the best dimensional constant in map
stability is at least of order $d^{1/3}$.
\end{abstract}

\section{Introduction and main results}
\label{sec:introduction}

Let $\mathcal P_2(\R^d)$ denote the probability measures on $\R^d$,
$d\geqslant2$, with finite second moment, equipped with the quadratic
Wasserstein distance $W_2$.
For $\rho=f\dd x\in\mathcal P_2(\R^d)$, write
$T_\mu=\nabla\phi_\mu$ for the quadratic optimal transport map from
$\rho$ to $\mu\in\mathcal P_2(\R^d)$, with $\phi_\mu$
convex~\cite{BrenierPolar,McCannMonotone}. Quantitative stability asks how
the potential, modulo additive constants, and the map change when the
target changes. We refer to these as \emph{potential stability} and
\emph{map stability}, respectively, and study them for log-concave
sources and bounded targets.

For source densities bounded above and below by positive constants
on a convex body, Delalande and Mérigot
\cite{DelalandeMerigotQuantitative} obtained stability of potentials
and converted it into stability of maps by a reverse Poincaré
inequality.
Letrouit and Mérigot~\cite{LetrouitMerigotGluing} extended stability
estimates to more general sources, including log-concave densities
with full support.

Mérigot~\cite[Theorem~3.1]{MerigotSharpStability} proved the sharp
$W_1^{1/4}$ map estimate for bounded, compactly supported source densities.
Letrouit~\cite{LetrouitUnstable} constructed uniform sources on nonconvex
domains for which no $W_2$ map estimate can have exponent larger than $1/3$.
For uniform sources on convex bodies, the sharp $W_2^{1/3}$ map estimate
was established in the OpenAI preprint~\cite{OpenAISharpOneThird}.

For general BV sources, we previously obtained in~\cite{sotr} the
$W_2^{1/4}$ map estimate under a support bound on one target, with
constants depending on the total variation of the source density.

Letrouit and Manole~\cite{LetrouitManoleBrenierStability} proved
linear $W_2$ stability of potentials modulo constants for strongly
log-concave sources and for log-concave sources on convex bodies,
using linearization of the Monge--Ampère equation. They call this
a \emph{Brenier stability inequality} and show that such an unweighted
bound can fail for general log-concave sources.
For log-concave sources on convex bodies, they also proved the sharp
$W_2^{1/3}$ map estimate when the negative log-density is finite
throughout the body~\cite[Corollary~1.5]{LetrouitManoleBrenierStability}.
Ruan~\cite{RuanPotentialMapStability}
proved linear potential stability for full-support strongly log-concave
sources by semi-discrete transport, and combined it with an unweighted
BV reverse Poincaré inequality to obtain the sharp map exponent $1/3$.
His interpolation inequality applies to arbitrary BV source densities
\cite[Appendix~A, Theorem~A.2]{RuanPotentialMapStability}.

We extend the $1/3$ map estimate to arbitrary log-concave sources.
We measure potential differences in a weighted $L^2$ norm that allows
them to grow in the tails of the source. For potential stability, we
adapt Letrouit and Manole's linearized Monge--Ampère argument to this
norm. For the reverse Poincaré inequality, we adapt Ruan's BV interpolation
argument to the same weight. Combining these weighted estimates gives
map stability, with a constant that grows polynomially in dimension.

\paragraph{AI use}
The approach developed in this paper emerged from discussions with
ChatGPT (6.1 Sol). Motivated by the OpenAI preprint~\cite{OpenAISharpOneThird},
the authors began these discussions on October 19, 2026, by asking whether
the $1/3$ stability estimate extends to Gaussian and strongly log-concave
sources. The resulting proof combined an unweighted potential stability
estimate as in Letrouit and Manole~\cite{LetrouitManoleBrenierStability}
with an unweighted BV interpolation inequality as in
Ruan~\cite[Appendix~A, Theorem~A.2]{RuanPotentialMapStability}.
Further discussions on October 22, aimed at extending the estimate to
general log-concave sources, led to the weighted estimates of this paper.
The authors checked and revised the manuscript and take full
responsibility for its content.

\paragraph{Funding}
This work benefited from support from the French State,
managed by the French National Research Agency (ANR)
as part of the France 2030 programme, reference ANR-23-IACL-0008.

\subsection{Weighted potential stability}
\label{subsec:weighted-brenier-stability}

We write the log-concave source as
\begin{equation}\label{eq:source}
 \rho=f\dd x,\qquad f=Z^{-1}e^{-V},\qquad \min V=0,\qquad
 Z=\int e^{-V}\dd x\in(0,\infty),
\end{equation}
where $V\colon\R^d\to[0,\infty]$ is lower semicontinuous and convex.
We call $\rho$ isotropic if it has mean zero and covariance $I_d$,
and write $T_\#\rho$ for its pushforward by $T$.

For a probability measure $\sigma$ and a measurable weight $U>0$
that is finite $\sigma$-almost everywhere, define the weighted variance
\begin{equation}\label{eq:weighted-variance}
 \Var_{\sigma}(w;U)
 \coloneqq\inf_{a\in\R}\int\frac{\lvert w-a\rvert^2}{U^2}\dd\sigma,
\end{equation}
provided $0<\int U^{-2}\dd\sigma<\infty$ and $w/U\in L^2(\sigma)$.
We use $U=\lambda+V$, $\lambda>0$,
with all integrals taken on $\{V<\infty\}$.
Thus $\Var_{\rho}(\phi-\psi;\lambda+V)^{1/2}$ is the weighted
$L^2(\rho)$ distance modulo additive constants.

In this norm, we have the following linear potential stability bound.

\begin{theorem}[Weighted potential stability]
\label{thm:weighted-stability}
Let $\rho=f\dd x$ be isotropic and log-concave on $\R^d$, with $V$
normalized as in \eqref{eq:source}. For any convex Lipschitz $\phi$, $\psi$ on
$\R^d$,
\[
 \Var_{\rho}(\phi-\psi;4d+V)^{1/2}
 \le3\sqrt e\,W_2\bigl((\nabla\phi)_\#\rho,(\nabla\psi)_\#\rho\bigr).
\]
\end{theorem}

We prove Theorem~\ref{thm:weighted-stability} in
Section~\ref{sec:weighted-stability}.

The linear growth in $V$ of the weight in
Theorem~\ref{thm:weighted-stability} is sharp within the family of
power weights.

\begin{proposition}[Sharpness of the weight exponent]
\label{prop:weight-power}
Let $\rho=e^{-V}\dd x$ on $\R^2$, where
\[
 V(x)=\begin{cases}
       \sqrt2x_2+2,&x_2+\sqrt2\ge\lvert x_1\rvert,\\
       +\infty,&\text{otherwise}.
      \end{cases}
\]
Then $\rho$ is isotropic and log-concave, with $\min V=0$.
For every $\lambda>0$ and $p<2$, there exist two-point targets
$\mu_n$, $\nu_n$ in $\overline B(0,1)$ such that
$W_2(\mu_n,\nu_n)\to0$ and, as $n\to\infty$,
\[
 \frac{\displaystyle\inf_{a\in\R}
       \int\frac{\lvert \phi_{\mu_n}-\phi_{\nu_n}-a\rvert^2}{(\lambda+V)^p}\dd\rho}
      {W_2(\mu_n,\nu_n)^2}
 \sim\frac16 n^{2-p}.
\]
\end{proposition}

We prove Proposition~\ref{prop:weight-power} in
Section~\ref{sec:weight-optimality}.

\paragraph{A Poincaré interpretation}
\label{subsec:poincare-interpretation}

Set $U=4d+V$. For any measurable $S$ with $S_\#\rho=\rho$,
the pair $(\nabla\phi,\nabla\psi\circ S)$ gives a coupling of
the gradient laws. Theorem~\ref{thm:weighted-stability} therefore gives
\[
 \Var_\rho(\phi-\psi;U)
 \le9e\int{\lvert \nabla\phi-\nabla\psi\circ S\rvert}^2\dd\rho.
\]
Taking $S=\mathrm{id}$ gives a weighted Poincaré bound for
$w=\phi-\psi$.

\subsection{Weighted reverse Poincaré inequality}
\label{subsec:weighted-interpolation}

To convert control of potentials into control of maps, we prove a
weighted reverse Poincaré inequality with the same weight.
We write $\BV(\R^d)$ for the integrable functions whose distributional
derivative $Dh$ is a finite vector measure, with total variation measure $\lvert Dh\rvert$.
Densities are extended by zero outside their support, so their total
variation includes boundary contributions.

\begin{theorem}[Weighted reverse Poincaré inequality]
\label{thm:weighted-interpolation}
Let $\rho=f\dd x$ be log-concave on $\R^d$, with $V$ normalized as in
\eqref{eq:source}. Then $f\in\BV(\R^d)$, and for any convex
$R$-Lipschitz $\phi$, $\psi$ on $\R^d$ and any $\lambda\ge1/2$,
\[
 \norm{\nabla\phi-\nabla\psi}_{L^2(\rho)}^3
 \le8\sqrt3\,(\lambda+d)\lvert Df\rvert(\R^d)R^2
       \Var_{\rho}(\phi-\psi;\lambda+V)^{1/2}.
\]
\end{theorem}

We prove Theorem~\ref{thm:weighted-interpolation} in
Section~\ref{sec:interpolation}.

\subsection{Map stability}

To express the map estimate in terms of dimension alone, we prove
an explicit version of the BV bound of Eldan and
Klartag~\cite[Lemma~11]{EldanKlartagPointwise}, using the sublevel sets
of $V$.

\begin{proposition}[Isotropic BV bound]
\label{prop:isotropic-bv}
Let $\rho=f\dd x$ be isotropic and log-concave on $\R^d$. Then
\[
 \lvert Df\rvert(\R^d)\le12d.
\]
\end{proposition}

We prove Proposition~\ref{prop:isotropic-bv} in
Appendix~\ref{sec:explicit-bv}.

\begin{corollary}[Map stability]
\label{cor:map-stability}
Let $\rho=f\dd x$ be isotropic and log-concave on $\R^d$, and let
$\mu$, $\nu$ be probability measures supported in a common ball of
radius $R>0$. Their Brenier maps satisfy
\[
 \begin{aligned}
 \norm{T_\mu-T_\nu}_{L^2(\rho)}
 &\le\bigl(120\sqrt{3e}\,d\lvert Df\rvert(\R^d)\bigr)^{\!1/3}R^{2/3}W_2(\mu,\nu)^{1/3}\\
 &\le(1440\sqrt{3e})^{1/3}d^{2/3}R^{2/3}W_2(\mu,\nu)^{1/3}.
 \end{aligned}
\]
\end{corollary}

\begin{proof}
Let $b$ be the centre of the common target ball. Choose globally
$R$-Lipschitz convex potentials for the translated targets, whose
gradients are $T_\mu-b$ and $T_\nu-b$. Applying
Theorems~\ref{thm:weighted-stability} and~\ref{thm:weighted-interpolation}
with $\lambda=4d$ gives the first inequality; Proposition~\ref{prop:isotropic-bv}
gives the second.
\end{proof}

Since $(1440\sqrt{3e})^{1/3}<17$, this yields the estimate stated in
the abstract.

\begin{remark}[General covariance]
\label{rem:covariance}
If $\rho$ is log-concave with mean $m$ and covariance $\Sigma\succ0$,
then
\[
 \norm{T_\mu-T_\nu}_{L^2(\rho)}
 \le (1440\sqrt{3e})^{1/3}d^{2/3}
 \sqrt{\frac{\lambda_{\max}(\Sigma)}{\lambda_{\min}(\Sigma)}}
 R^{2/3}W_2(\mu,\nu)^{1/3}.
\]
This follows from Corollary~\ref{cor:map-stability} after the change of variables
$x=m+\Sigma^{1/2}z$, with targets $(\Sigma^{1/2})_\#\mu$ and
$(\Sigma^{1/2})_\#\nu$.
\end{remark}

\subsection{Dimension dependence}

Let $C_d^{\mathrm{best}}$ be the least constant such that
\[
 \norm{T_\mu-T_\nu}_{L^2(\rho)}
 \le C_d^{\mathrm{best}}W_2(\mu,\nu)^{1/3}
\]
for every isotropic log-concave source on $\R^d$ and all targets
in $\overline B(0,1)$.

\begin{proposition}[Sharpness and dimensional lower bound]
\label{prop:dimension-lower-bound}
For every $d\ge2$, set $k_d=2\lfloor d/2\rfloor$ and
\[
 \rho=\operatorname{Unif}([-\sqrt3,\sqrt3]^d).
\]
There exist finitely supported targets $\mu_n$, $\nu_n$ in
$\overline B(0,1)$ whose Brenier maps from $\rho$ satisfy, as $n\to\infty$,
\[
 W_2(\mu_n,\nu_n)^2\sim\frac{k_d}{4n^3},\qquad
 \norm{T_{\mu_n}-T_{\nu_n}}_{L^2(\rho)}^2
 \ge(1-o(1))\frac{k_d}{4n}.
\]
\end{proposition}

We prove Proposition~\ref{prop:dimension-lower-bound} in
Section~\ref{sec:dimension-lower-bounds}. Its estimates show that the
exponent $1/3$ in Corollary~\ref{cor:map-stability} is optimal even for
this fixed source and give $C_d^{\mathrm{best}}\ge(k_d/4)^{1/3}$.
Since $k_d\ge2d/3$, we obtain, together with
Corollary~\ref{cor:map-stability},
\[
 6^{-1/3}d^{1/3}\le C_d^{\mathrm{best}}
 \le(1440\sqrt{3e})^{1/3}d^{2/3},\qquad d\ge2.
\]
The optimal power of the dimension between these bounds is open.

\section{Weighted potential stability}
\label{sec:weighted-stability}

We establish a weighted coercivity estimate for the linearized
Monge--Ampère equation, adapting the unweighted argument of Letrouit
and Manole~\cite{LetrouitManoleBrenierStability}. The weight depends on
the source potential $V$, for which we first prove a geometric estimate.

\begin{lemma}[Source geometry]
\label{lem:source-geometry}
Let $\rho$ be isotropic and log-concave, with $V$ normalized as in
\eqref{eq:source}. Then
\[
 \lvert x\rvert\le\sqrt e\,\bigl(2d+V(x)\bigr)
 \qquad\text{for }x\in\{V<\infty\}.
\]
If $\phi$, $\psi\colon\R^d\to\R$ are $R$-Lipschitz and satisfy
$\phi(0)=\psi(0)=0$, then
\[
 \frac{\lvert \phi-\psi\rvert}{4d+V}\le2\sqrt e\,R
 \qquad\text{on }\{V<\infty\}.
\]
\end{lemma}

\begin{proof}
Fix $x\in\{V<\infty\}$, put $\ell=V(x)$ and
$\theta=(2d+\ell)^{-1}\le1/2$, and define
$A_\theta(y)=(1-\theta)y+\theta x$.
Convexity and $V\ge0$ give
$f(A_\theta(y))\ge e^{-\theta\ell}f(y)$.
A change of variables therefore shows that
\[
 \rho_\theta\coloneqq(A_\theta)_\#\rho
 \le(1-\theta)^{-d}e^{\theta\ell}\rho\le e\rho,
\]
since $-\log(1-\theta)\le2\theta$ and
$(2d+\ell)\theta=1$. If $x\ne0$, set $n=x/\lvert x\rvert$.
The mean of $\rho_\theta$ is $\theta x$, so Jensen's inequality and
isotropy give
\[
 \theta^2\lvert x\rvert^2
 \le\int(n\cdot z)^2\dd\rho_\theta(z)
 \le e\int(n\cdot z)^2\dd\rho(z)=e.
\]
This proves the radius bound.
Finally, $\lvert \phi(x)-\psi(x)\rvert\le2R\lvert x\rvert$ proves the weighted bound.
\end{proof}

The following weighted coercivity estimate is the main technical
contribution of this paper.

\begin{lemma}[Weighted divergence estimate]
\label{lem:weighted-divergence}
Let $\Omega\subset\R^d$ be a bounded smooth convex domain, let $V\ge0$
be smooth and convex on $\overline\Omega$, and put
$Z_\Omega=\int_\Omega e^{-V}\dd x$ and
$\rho_\Omega=Z_\Omega^{-1}e^{-V}\dd x$ on $\Omega$.
Suppose that $U=4d+V$ satisfies $\lvert x\rvert\le LU$ on $\Omega$ for some $L>0$.
Let $\phi$ be smooth on $\overline\Omega$, with $H=D^2\phi\succ0$,
and let $\chi\in\mathcal C^3(\overline\Omega)$.
Set
\[
 T=H^{-1}\nabla \chi,\qquad
 q=-e^V\diver(e^{-V}T).
\]
If $T\cdot n=0$ on $\partial\Omega$, where $n$ is the outward unit
normal, then
\[
 \norm{T}_{L^2(\rho_\Omega)}
 \le3L\norm{Uq}_{L^2(\rho_\Omega)}.
\]
\end{lemma}

\begin{proof}
Put $\tau=\tr((DT)^2)$, with $(DT)_{ij}=\partial_jT_i$.
Differentiation of $HT=\nabla \chi$ gives
\[
 HDT=D^2\chi-D^3\phi[T],
\]
where $(D^3\phi[T])_{ij}=\sum_\ell\partial_{ij\ell}\phi\,T_\ell$.
Thus $H^{1/2}DTH^{-1/2}$ is symmetric. Its eigenvalues are real,
so Cauchy--Schwarz gives
\begin{equation}\label{eq:div-trace-control}
 \tau\ge0,\qquad (\diver T)^2\le d\tau.
\end{equation}

Write $\diver_{\rho_\Omega}G=\diver G-\nabla V\cdot G$.
Direct differentiation gives
\[
 T\cdot\nabla q
 =-\diver_{\rho_\Omega}((DT)T)+\tau+D^2V[T,T].
\]
For a smooth scalar function $\Phi$, substitute this into
$\int\Phi q^2\dd\rho_\Omega
 =\int T\cdot\nabla(\Phi q)\dd\rho_\Omega$ and integrate by parts.
Using
\[
 \nabla\Phi\cdot(DT)T
 =T\cdot\nabla(\nabla\Phi\cdot T)-D^2\Phi[T,T]
\]
and $((DT)T)\cdot n=-\mathrm{II}(T,T)$ on $\partial\Omega$, we obtain
\begin{equation}
\begin{multlined}
 \int_\Omega\Phi q^2\dd\rho_\Omega
 =\int_\Omega\bigl\{
       \Phi\tau+(\Phi D^2V-D^2\Phi)[T,T]
       +2q\nabla\Phi\cdot T\bigr\}\dd\rho_\Omega \\
 +\int_{\partial\Omega}\Phi\,\mathrm{II}(T,T)
          \frac{e^{-V}}{Z_\Omega}\dd S.
\end{multlined}
\label{eq:div-weighted-identity}
\end{equation}
Here $M[a,a]=a^{\mathsf T}Ma$ and
$\mathrm{II}(a,a)=a\cdot D_an$ for tangent vectors $a$.
The boundary identity follows by differentiating $T\cdot n=0$
in the tangent direction $T$. Convexity of $\Omega$ makes the
boundary term nonnegative whenever $\Phi\ge0$.

We now take $c=2L^2$ and
\[
 \Phi=cU^2-\lvert x\rvert^2.
\]
Since $U\ge4$, $U^2\ge8d$ and $\lvert x\rvert\le LU$,
\[
 \Phi\ge\frac c2 U^2\ge4cd,
 \qquad \Phi-2cU\ge0.
\]
Write $b=\diver T$ and $s=\nabla V\cdot T=q+b$.
The volume integrand in \eqref{eq:div-weighted-identity} becomes
\begin{align*}
 &\Phi\tau+(\Phi-2cU)D^2V[T,T]+2\lvert T\rvert^2\\
 &\hspace{35mm}-2cs^2+4cUqs-4q\,x\cdot T.
\end{align*}
The Hessian and boundary terms are nonnegative. Moreover,
the trace bound \eqref{eq:div-trace-control} and $\Phi\ge4cd$ give
$\Phi\tau\ge4cb^2$. Completing the square yields
\begin{align*}
 \Phi\tau-2cs^2+4cUqs
 &\ge2c\bigl(b+(U-1)q\bigr)^2
       +2c(-U^2+4U-2)q^2\\
 &\ge-2cU^2q^2.
\end{align*}
Discarding the nonnegative terms in the weighted identity gives
\[
 2\int_\Omega{\lvert T\rvert}^2\dd\rho_\Omega
 \le\int_\Omega\bigl((\Phi+2cU^2)q^2+4q\,x\cdot T\bigr)
       \dd\rho_\Omega.
\]
Finally use $4\lvert q\,x\cdot T\rvert\le\lvert T\rvert^2+4\lvert x\rvert^2q^2$ to obtain
\[
 \int_\Omega{\lvert T\rvert}^2\dd\rho_\Omega
 \le\int_\Omega(3cU^2+3{\lvert x\rvert}^2)q^2\dd\rho_\Omega
 \le9L^2\int_\Omega U^2q^2\dd\rho_\Omega.
\]
Taking square roots proves the estimate.
\end{proof}

We apply the divergence estimate by duality to control the variation
of a Brenier potential along a smooth target path.

\begin{lemma}[Variation of potentials along a smooth path]
\label{lem:path-potential-variation}
Let $\Omega$ and $\mathcal D$ be balls centred at the origin, with
outward unit normals $n_\Omega$, $n_{\mathcal D}$, and let
$\rho_\Omega=Z_\Omega^{-1}e^{-V}\dd x$, where $V\ge0$ is smooth and
convex on $\overline\Omega$. Let $h_t>0$ be a smooth density and
$j_t$ a smooth $\R^d$-valued flux on
$[0,1]\times\overline{\mathcal D}$ satisfying
\[
 \int_{\mathcal D}h_t\dd y=1,\qquad
 \partial_t h_t+\diver j_t=0,
 \qquad j_t\cdot n_{\mathcal D}=0\quad\text{on }\partial\mathcal D.
\]
Set $\mu_t=h_t\dd y$, define the velocity $v_t=j_t/h_t$, and let $\phi_t$ be the
Brenier potential from $\rho_\Omega$ to $\mu_t$.
If $U=4d+V$ satisfies $\lvert x\rvert\le\sqrt e\,U$ on $\Omega$, then
\[
 \Var_{\rho_\Omega}(\phi_1-\phi_0;U)^{1/2}
 \le3\sqrt e\int_0^1\norm{v_t}_{L^2(\mu_t)}\dd t.
\]
\end{lemma}

\begin{proof}
The regularity theory for the second boundary value problem gives
smooth strictly convex potentials, with
$\nabla\phi_t\colon\overline\Omega\to\overline{\mathcal D}$ a
diffeomorphism; see \cite[Theorems~1.2--1.3]{TrudingerWangSecondBoundary}
and the smooth-data regularity following those theorems.
By \cite[Corollary~3.5]{GonzalezSanzShengLinearization},
the potentials form a $\mathcal C^1$ curve in
$\mathcal C^2(\overline\Omega)$ after normalization by additive constants.
Write $H_t=D^2\phi_t$ and $f_\Omega=Z_\Omega^{-1}e^{-V}$.

The pushforward relation gives the logarithmic Monge--Ampère equation
\[
 \log\det D^2\phi_t+\log h_t(\nabla\phi_t)=\log f_\Omega.
\]
At each fixed $t$, its linearization in the potential, with $h_t$ fixed,
is the operator
\[
 \begin{aligned}
 \mathcal L_t u
 &\coloneqq\tr(H_t^{-1}D^2u)
   +\nabla\log h_t(\nabla\phi_t)\cdot\nabla u\\
 &=f_\Omega^{-1}\diver(f_\Omega H_t^{-1}\nabla u).
 \end{aligned}
\]
The divergence form follows from the divergence-free cofactor identity
and $f_\Omega=h_t(\nabla\phi_t)\det H_t$.

We use this operator in a duality argument to estimate
$\dot\phi_t=\partial_t\phi_t$ in the weighted norm.
Normalize the potentials by $\int_\Omega \phi_t/U^2\dd\rho_\Omega=0$
and, at a fixed time $t$, put $g=\dot\phi_t/U^2$.
Then $\int g\dd\rho_\Omega=0$, and the uniformly
elliptic Neumann problem
\[
 -\mathcal L_t\chi=g
 \quad\text{in }\Omega,
 \qquad H_t^{-1}\nabla \chi\cdot n_\Omega=0
 \quad\text{on }\partial\Omega
\]
has a $\mathcal C^3$ solution unique modulo constants. Set $T=H_t^{-1}\nabla\chi$,
so that $q=g$ in Lemma~\ref{lem:weighted-divergence}.
When differentiating the pushforward relation at this time, take the
fixed test function $\zeta=\chi\circ(\nabla\phi_t)^{-1}$. Its gradient satisfies
$\nabla\zeta(\nabla\phi_t)=T$. Integration by parts and the continuity
equation give
\[
 \norm{\dot\phi_t/U}_{L^2(\rho_\Omega)}^2
 =\int_\Omega g\dot\phi_t\dd\rho_\Omega
   =\int_\Omega T\cdot\nabla\dot\phi_t\dd\rho_\Omega
 =\int_\Omega T\cdot v_t(\nabla\phi_t)\dd\rho_\Omega.
\]
By Lemma~\ref{lem:weighted-divergence} with $L=\sqrt e$,
$\norm{T}_{L^2(\rho_\Omega)}\le3\sqrt e\,\norm{Ug}_{L^2(\rho_\Omega)}
=3\sqrt e\,\norm{\dot\phi_t/U}_{L^2(\rho_\Omega)}$.
Applying Cauchy--Schwarz to the preceding identity yields
\[
 \norm{\dot\phi_t/U}_{L^2(\rho_\Omega)}
 \le3\sqrt e\,\norm{v_t}_{L^2(\mu_t)}.
\]
Integrating in $t$ and applying the triangle inequality in
$L^2(\rho_\Omega)$ proves the claimed bound.
\end{proof}

A Wasserstein geodesic has length $W_2(\mu,\nu)$ but may not satisfy
the smoothness and positivity assumptions of
Lemma~\ref{lem:path-potential-variation}. In the following proof, we
regularize it without increasing its length.

\begin{lemma}[Weighted potential stability on a source ball]
\label{lem:bounded-source-stability}
Let $\Omega$ be a ball centred at the origin, and let
$\rho_\Omega=Z_\Omega^{-1}e^{-V}\dd x$, where $V\ge0$ is smooth and
convex on $\overline\Omega$. Suppose that $U=4d+V$ satisfies
$\lvert x\rvert\le\sqrt e\,U$ on $\Omega$. Then, for any convex Lipschitz
$\phi$, $\psi$ on $\R^d$,
\[
 \Var_{\rho_\Omega}(\phi-\psi;U)^{1/2}
 \le3\sqrt e\,W_2\bigl((\nabla\phi)_\#\rho_\Omega,
                       (\nabla\psi)_\#\rho_\Omega\bigr).
\]
\end{lemma}

\begin{proof}
Subtract constants so that $\phi(0)=\psi(0)=0$, choose a common
Lipschitz bound $R>0$, and set
$\mu=(\nabla\phi)_\#\rho_\Omega$ and
$\nu=(\nabla\psi)_\#\rho_\Omega$.
Let $\pi$ be a coupling attaining their quadratic transport cost.
For $0<\varepsilon<1/2$, choose a nonnegative smooth probability kernel
$\kappa_\varepsilon$ supported in $B(0,\varepsilon)$, put
$\mathcal D=B(0,R+2\varepsilon)$, and write
$y_t=(1-t)y_0+ty_1$ for $(y_0,y_1)$ distributed according to $\pi$.
For $0<\alpha<1$, define on $\mathcal D$
\[
 \begin{aligned}
 h_t(y)&=(1-\alpha)\int
       \kappa_\varepsilon(y-y_t)\dd\pi(y_0,y_1)
       +\frac\alpha{\lvert \mathcal D\rvert},\\
 j_t(y)&=(1-\alpha)\int
       (y_1-y_0)\kappa_\varepsilon(y-y_t)\dd\pi(y_0,y_1),
 \qquad v_t=\frac{j_t}{h_t}.
 \end{aligned}
\]
These functions are smooth up to the boundary, $h_t\ge\alpha/\lvert \mathcal D\rvert$,
and $\int_{\mathcal D}h_t\dd y=1$. Since ${\lvert y_t\rvert}\le R$, the flux vanishes
on a neighbourhood of $\partial\mathcal D$. Differentiation under the
integral gives $\partial_t h_t+\diver j_t=0$.
Pointwise Cauchy--Schwarz gives
\[
 \frac{\lvert j_t(y)\rvert^2}{h_t(y)}
 \le(1-\alpha)\int{\lvert y_1-y_0\rvert}^2
               \kappa_\varepsilon(y-y_t)\dd\pi(y_0,y_1).
\]
Integrating over $\mathcal D$ therefore yields
\[
 \int_{\mathcal D}{\lvert v_t\rvert}^2h_t\dd y
 \le(1-\alpha)\int{\lvert y_1-y_0\rvert}^2\dd\pi(y_0,y_1)
 \le W_2(\mu,\nu)^2.
\]
The path length is at most $W_2(\mu,\nu)$.
Its endpoints converge in $W_2$ to $\mu$, $\nu$ as
$\alpha,\varepsilon\downarrow0$: convolution and the vanishing mixture
weight give weak convergence, and all the endpoint measures are supported
in the fixed ball $\overline B(0,R+1)$.

We can now apply Lemma~\ref{lem:path-potential-variation} to the regularized
path $h_t\dd y$. Denote its potentials by
$\phi_t^{\varepsilon,\alpha}$. The length bound gives
\[
 \Var_{\rho_\Omega}
   (\phi_1^{\varepsilon,\alpha}-\phi_0^{\varepsilon,\alpha};U)^{1/2}
 \le3\sqrt e\,W_2(\mu,\nu).
\]
Normalize these potentials at zero. Their gradients lie in
$\overline B(0,R+2\varepsilon)$, so the potentials are uniformly
Lipschitz on $\overline\Omega$. Arzelà--Ascoli gives a subsequence
converging uniformly there as $\varepsilon,\alpha\downarrow0$.
For convex functions this implies almost everywhere convergence of
the gradients in $\Omega$. Dominated convergence identifies their
limit laws as $\mu$ and $\nu$. Brenier uniqueness and positivity of
the source density identify the limiting potentials with $\phi$
and $\psi$ on the connected source ball, respectively; normalization
at zero fixes the additive constants.

For any weighted variance as in \eqref{eq:weighted-variance},
\begin{equation}\label{eq:path-weighted-variance}
 \Var_{\sigma}(w;U)
 =\int\frac{w^2}{U^2}\dd\sigma
  -\frac{\bigl(\int w/U^2\dd\sigma\bigr)^2}
         {\int U^{-2}\dd\sigma}.
\end{equation}
Uniform convergence on $\overline\Omega$ and $U\ge4$ pass all three
integrals to the limit, proving the claim.
\end{proof}

The following source approximation controls both the density and its
normalized source potential.

\begin{lemma}[Smooth isotropic approximation]
\label{lem:source-approximation}
Let $\rho=f\dd x$ be isotropic and log-concave. For independent
$X\sim\rho$ and $G\sim N(0,I_d)$, define
\[
 a_\eta=\sqrt{1+\eta^2},\qquad
 X_\eta=\frac{X+\eta G}{a_\eta},\qquad
 \rho_\eta=f_\eta\dd x,\qquad \eta>0.
\]
Then $f_\eta$ is smooth, positive, isotropic and log-concave, and
\[
 f_\eta\longrightarrow f\quad\text{in }L^1(\R^d),\qquad
 \norm{f_\eta}_\infty\longrightarrow\norm f_\infty.
\]
If $V_\eta=\log{\norm{f_\eta}_\infty}-\log f_\eta$, then
$V_\eta\to V$ at every interior point of $\supp\rho$.
\end{lemma}

\begin{proof}
Let $\gamma_\eta$ be the density of $\eta G$. Then
\[
 f_\eta(x)=a_\eta^d(f*\gamma_\eta)(a_\eta x).
\]
Smoothness and positivity follow from Gaussian convolution.
Prékopa's theorem gives log-concavity. The normalization gives
\[
 \E X_\eta=0,\qquad
 \Cov(X_\eta)=a_\eta^{-2}(I_d+\eta^2I_d)=I_d.
\]
The approximate-identity property and continuity of dilations on $L^1$
give
\[
 \norm{f_\eta-f}_{L^1}
 \le\norm{f*\gamma_\eta-f}_{L^1}
   +\norm{a_\eta^df(a_\eta\,\cdot)-f}_{L^1}\longrightarrow0.
\]
A log-concave density is bounded and continuous in the interior of
its support. Locally uniform approximate-identity convergence there gives
$f_\eta(x)\to f(x)$. Since
$\norm{f_\eta}_\infty\le a_\eta^d\norm f_\infty$, and interior values
of $f$ approach its supremum, the suprema converge as well.
This proves the assertion about $V_\eta$.
\end{proof}

\begin{proof}[Proof of Theorem~\ref{thm:weighted-stability}]
Subtract constants so that $\phi(0)=\psi(0)=0$, and let $R$ be a common
Lipschitz bound. Let $f_\eta$ and $V_\eta$ be given by
Lemma~\ref{lem:source-approximation}, and put
\[
 \begin{gathered}
 U_\eta=4d+V_\eta,\qquad
 \Omega_\eta=B(0,\eta^{-1}),\qquad
 p_\eta=\int_{\Omega_\eta}f_\eta\dd x,\\
 \widetilde\rho_\eta=\widetilde f_\eta\dd x
 \coloneqq p_\eta^{-1}\one_{\Omega_\eta}f_\eta\dd x.
 \end{gathered}
\]
Lemma~\ref{lem:source-geometry}, applied to the isotropic density
$f_\eta$, gives $\lvert x\rvert\le\sqrt e\,U_\eta$ on $\Omega_\eta$.
Thus Lemma~\ref{lem:bounded-source-stability}, with the restriction of
$V_\eta$ to $\Omega_\eta$, gives
\[
 \Var_{\widetilde\rho_\eta}(\phi-\psi;U_\eta)^{1/2}
 \le3\sqrt e\,W_2\bigl((\nabla\phi)_\#\widetilde\rho_\eta,
                       (\nabla\psi)_\#\widetilde\rho_\eta\bigr).
\]
Since $f_\eta\to f$ in $L^1(\R^d)$ and the balls $\Omega_\eta$ exhaust
$\R^d$, we have $p_\eta\to1$ and $\widetilde f_\eta\to f$ in $L^1$.
The gradient laws therefore converge in total variation to
$(\nabla\phi)_\#\rho$ and $(\nabla\psi)_\#\rho$.
All these laws lie in $\overline B(0,R)$, so convergence also holds in $W_2$.

Put $w=\phi-\psi$ and $U=4d+V$. The source geometry gives
\[
 \lvert w\rvert/U_\eta\le2\sqrt e\,R,\qquad U_\eta^{-1}\le1/4.
\]
Each of $w^2/U_\eta^2$, $w/U_\eta^2$ and $1/U_\eta^2$ is therefore
uniformly bounded, and converges pointwise in the interior of
$\supp\rho$. Its integral against $\widetilde f_\eta-f$ tends to zero by the
$L^1$ convergence of the densities; its integral against $f$ converges
by dominated convergence. The boundary of the convex support has
Lebesgue measure zero. Formula~\eqref{eq:path-weighted-variance},
whose limiting denominator is positive, now gives
\[
 \Var_{\widetilde\rho_\eta}(\phi-\psi;U_\eta)
 \longrightarrow\Var_{\rho}(\phi-\psi;U).
\]
Passing to the limit proves the theorem.
\end{proof}

\section{Weighted reverse Poincaré inequality}
\label{sec:interpolation}

We adapt the BV interpolation argument
of Ruan~\cite[Appendix~A, Theorem~A.2]{RuanPotentialMapStability}
by using a truncation threshold proportional to $\lambda+V$.
The density $g_\lambda=(\lambda+1+V)f$ satisfies
$\nabla g_\lambda=(\lambda+V)\nabla f$ in the smooth case.
We first bound its total variation, which controls the additional
integration-by-parts terms arising from the variable threshold.

\begin{lemma}[Weighted BV bound]
\label{lem:weighted-bv}
Let $\rho=f\dd x$ be log-concave, with $V$ normalized as in
\eqref{eq:source}. For every $\lambda>0$, the functions $f$ and
$g_\lambda=(\lambda+1+V)f$, extended by zero where $f=0$, belong
to $\BV(\R^d)\cap L^1(\R^d)$, and
\[
 \lvert Dg_\lambda\rvert(\R^d)\le(\lambda+d)\lvert Df\rvert(\R^d).
\]
\end{lemma}

\begin{proof}
Choose a minimizer $x_0$ of $V$, so that $V(x_0)=0$.
Starting from a sublevel set of positive volume,
contraction towards $x_0$ shows that every
$K_t=\{V\le t\}$, $t>0$, has nonempty interior. Since
$\lvert K_t\rvert\le e^tZ$, each $K_t$ is a convex body.
Writing $p(t)=\Per(K_t)=\lvert D\one_{K_t}\rvert(\R^d)$, convexity and
perimeter monotonicity give
\[
 x_0+\frac{s}{t}(K_t-x_0)\subset K_s,\qquad
 p(s)\le p(t)\le(t/s)^{d-1}p(s),\qquad 0<s<t.
\]
Thus $p$ is locally Lipschitz, $tp'(t)\le(d-1)p(t)$ almost
everywhere, and $p$ is bounded near zero and grows at most as
$t^{d-1}$ at infinity. Integration by parts yields
\[
 \int_0^\infty te^{-t}p(t)\dd t
 =\int_0^\infty e^{-t}\bigl(p(t)+tp'(t)\bigr)\dd t
 \le d\int_0^\infty e^{-t}p(t)\dd t.
\]
The contraction also bounds $\lvert K_t\rvert$ by $O(1+t^d)$, so the
layer-cake representations of $f$ and $g_\lambda$ are integrable.
Together with Gauss--Green, the perimeter bounds give finite
total variation. Coarea then gives
\[
 \lvert Df\rvert(\R^d)=\frac1Z\int_0^\infty e^{-t}p(t)\dd t,\qquad
 \lvert Dg_\lambda\rvert(\R^d)
 =\frac1Z\int_0^\infty(\lambda+t)e^{-t}p(t)\dd t,
\]
since $t\mapsto(\lambda+1+t)e^{-t}$ is strictly decreasing.
The preceding moment bound gives the claimed estimate on $\lvert Dg_\lambda\rvert(\R^d)$.
\end{proof}

\begin{lemma}[Weighted reverse Poincaré inequality]
\label{lem:bv-interpolation}
Let $f$ be a probability density with $0\le f\le M<\infty$.
For $\lambda\ge1/2$, set
\[
 U=\lambda+\log(M/f),\qquad g=(U+1)f,
\]
with $g=0$ where $f=0$. If $g\in\BV(\R^d)\cap L^1(\R^d)$,
then, for convex $R$-Lipschitz $\phi$, $\psi$ with
$(\phi-\psi)/U\in L^2(f\dd x)$,
\[
 \norm{\nabla\phi-\nabla\psi}_{L^2(f\dd x)}^3
 \le8\sqrt3\,R^2\lvert Dg\rvert(\R^d)
       \Var_{f\dd x}(\phi-\psi;U)^{1/2}.
\]
\end{lemma}

\begin{proof}
Put $A=\lvert Dg\rvert(\R^d)$ and $w=\phi-\psi-a$, with a fixed $a\in\R$.
First suppose that $f$ is smooth and positive and that the potentials
are smooth with bounded Hessians. For $\varepsilon>0$, set
\[
 C_\varepsilon=f\max\{-\varepsilon U,
                              \min\{w,\varepsilon U\}\}.
\]
Since $\nabla g=U\nabla f$, the chain rule gives
\[
 \nabla C_\varepsilon
 =f\one_{\{\lvert w\rvert<\varepsilon U\}}\nabla w+b\nabla g,
 \qquad
 b=\begin{cases}
       w/U,&\lvert w\rvert<\varepsilon U,\\
       \varepsilon\operatorname{sgn}(w)(U-1)/U,
           &\lvert w\rvert\ge\varepsilon U.
    \end{cases}
\]
Here $\lvert b\rvert\le\varepsilon$ because $U\ge1/2$, and
$\lvert C_\varepsilon\rvert\le\varepsilon g$. Convexity and integration by
parts therefore yield
\begin{align}
 \int_{\{\lvert w\rvert<\varepsilon U\}}\lvert \nabla w\rvert^2f\dd x
 &=-\int C_\varepsilon\Delta w\dd x
                    -\int b\nabla w\cdot\nabla g\dd x\notag\\
 &\le\varepsilon\int g(\Delta\phi+\Delta\psi)\dd x
                    +2R\varepsilon A
 \le4R\varepsilon A.\label{eq:weighted-truncated-energy}
\end{align}
In the last step,
$0\le\int g\Delta\phi=-\int\nabla\phi\cdot\nabla g\le RA$,
and similarly for $\psi$. Spatial cutoffs justify the integrations
on $\R^d$, since $g$, $C_\varepsilon\in L^1$ and the gradients are bounded.

To obtain \eqref{eq:weighted-truncated-energy} for a general $g$,
write $H(s)=s(\lambda+1+\log(M/s))$, with $H(0)=0$.
Then $H'(s)=\lambda+\log(M/s)>0$ on $(0,M]$.
Convolve $g$ with Gaussian approximate identities to obtain
$0<g_n\le H(M)$, with $\lvert Dg_n\rvert(\R^d)\le A$ and $g_n\to g$
almost everywhere along a subsequence. Set $f_n=H^{-1}(g_n)$.
Since $f_n\le g_n/(\lambda+1)$ is integrable, the preceding calculation
applies to $f_n$ without requiring $\int f_n\dd x=1$.
On $\{f>0\}$, $f_n\to f$ and
$\lambda+\log(M/f_n)\to U$. Fatou's lemma on the strict truncation
region proves \eqref{eq:weighted-truncated-energy}.

Let $\dd\xi=\lvert \nabla w\rvert^2f\dd x$, $m=\xi(\R^d)$,
$L=4RA$ and $r=\lvert w\rvert/U$. Then
\eqref{eq:weighted-truncated-energy} says $\xi\{r<t\}\le Lt$.
If $L=0$, letting $t\to\infty$ gives $m=0$.
Otherwise, since $\dd\xi\le4R^2f\dd x$, the layer-cake formula gives
\[
 \int r^2f\dd x
 \ge\frac1{4R^2}\int r^2\dd\xi
 \ge\frac1{4R^2}\int_0^{m/L}2t(m-Lt)\dd t
 =\frac{m^3}{12R^2L^2}.
\]
Taking square roots and then the infimum over $a$ proves the claimed
inequality for smooth potentials. To extend it to general
convex Lipschitz potentials, we convolve both with a smooth compactly
supported approximate identity.
The Lipschitz bounds are preserved, the potentials converge uniformly,
and the gradients converge almost everywhere. Since $U\ge\lambda$,
the final inequality passes to the limit by dominated convergence.
\end{proof}

\begin{proof}[Proof of Theorem~\ref{thm:weighted-interpolation}]
With $M=\norm f_\infty=Z^{-1}$, the weight in
Lemma~\ref{lem:bv-interpolation} is $U=\lambda+V$.
Its integrability assumption holds because the potentials are
Lipschitz, $\lambda\ge1/2$ and $\rho\in\mathcal P_2(\R^d)$.
Combining Lemmas~\ref{lem:weighted-bv} and~\ref{lem:bv-interpolation}
proves the theorem.
\end{proof}

\section{Sharpness of the weight exponent}
\label{sec:weight-optimality}

\begin{proof}[Proof of Proposition~\ref{prop:weight-power}]
We first work before isotropic normalization. Let $(Y_0,S)$ have
law $\widehat\rho$ with log-concave density
\[
 \frac12 e^{-s}\one_{\{s\ge\lvert z\rvert\}}.
\]
The marginal density of $S$ is $se^{-s}\one_{\{s\ge0\}}$;
conditionally on $S=s$, $Y_0$ is uniform on $[-s,s]$.
Thus the mean is $(0,2)$ and the covariance is $2I_2$.

Fix $t>0$, put $m_t=(t+1)e^{-t}$ and
$n_\theta=(\sin\theta,\cos\theta)$, and choose $c_\theta$ so that
\[
 \Prob(n_\theta\cdot(Y_0,S)>c_\theta)=m_t.
\]
The half-space mass is smooth near $(\theta,c)=(0,t)$.
Its $c$ derivative there is $-te^{-t}$ and its $\theta$ derivative
vanishes by reflection in $z$. The implicit function theorem gives
$c_0=t$ and $c'_0=0$.
The convex $1$-Lipschitz functions
$\widehat\phi_\theta(z,s)=(n_\theta\cdot(z,s)-c_\theta)_+$ have gradient laws
\[
 \mu_\theta=(1-m_t)\delta_0+m_t\delta_{n_\theta},\qquad
 W_2(\mu_\theta,\mu_0)^2=m_t\lvert n_\theta-n_0\rvert^2
 \sim m_t\theta^2.
\]
For small $\theta$, matching the two corresponding atoms is optimal:
exchanging a mass $\alpha$ through zero adds $2\alpha\cos\theta$ to
its cost.

For each fixed $t$, as $\theta\to0$, the difference quotient converges in
$L^2((\lambda+s)^{-p}\widehat\rho)$ to $z\one_{\{s>t\}}$.
Indeed, it converges pointwise outside $\{s=t\}$ and is bounded by
$\lvert (z,s)\rvert+C_t$; all polynomial moments under this weighted measure
are finite. The limit has weighted mean zero by reflection.
Minimizing over constants is orthogonal projection in this $L^2$
space, so
\[
 \lim_{\theta\to0}
 \frac{\displaystyle\inf_a\int
       {\lvert \widehat\phi_\theta-\widehat\phi_0-a\rvert}^2(\lambda+s)^{-p}\dd\widehat\rho}
      {W_2(\mu_\theta,\mu_0)^2}
 =\frac{\displaystyle\int_t^\infty s^3e^{-s}(\lambda+s)^{-p}\dd s}
        {3(t+1)e^{-t}}.
\]
The factor $s^3/3$ comes from
$\frac12\int_{-s}^s z^2\dd z$. We next let $t\to\infty$.
Substituting $s=t+r$ and applying dominated convergence shows that
the right-hand side is asymptotic to $\frac13 t^{2-p}$.

Now set $x=(z,s-2)/\sqrt2$. The transformed source $\rho$ is isotropic,
with density $e^{-V}$ for the potential $V$ in the statement.
In particular, $\min V=0$. The potentials
$\phi_{\mu_\theta}(x)=\widehat\phi_\theta(\sqrt2x_1,\sqrt2x_2+2)/\sqrt2$
have the same gradient laws. For this source, the limit of the ratio as
$\theta\to0$ is half the expression above, hence is asymptotic to
$\frac16t^{2-p}$ as $t\to\infty$.

For each fixed $\lambda>0$ and $p<2$, take $t=n$ and choose
$0<\theta_n<1/n$ so that the ratio for the isotropic source differs
from its $\theta\to0$ limit by a relative error at most $1/n$.
With $m_n=(n+1)e^{-n}$, set
\[
 \mu_n=(1-m_n)\delta_0+m_n\delta_{n_{\theta_n}},\qquad
 \nu_n=(1-m_n)\delta_0+m_n\delta_{n_0}.
\]
The preceding asymptotic gives the claimed ratio, while
$W_2(\mu_n,\nu_n)^2=m_n\lvert n_{\theta_n}-n_0\rvert^2\to0$.
\end{proof}

\section{Dimensional lower bound}
\label{sec:dimension-lower-bounds}

\begin{proof}[Proof of Proposition~\ref{prop:dimension-lower-bound}]
First let $d$ be even. Take $0<\varepsilon<1/2$, write
$c=\sqrt3(1-\varepsilon)$, $\theta=\varepsilon/2$, and let $Q$ rotate
each coordinate plane $\operatorname{span}(e_{2j-1},e_{2j})$ through
$\theta$. Define convex $1$-Lipschitz potentials
\[
 \phi_\mu(x)=\max\{c,\lvert x_1\rvert,\ldots,\lvert x_d\rvert\},\qquad \phi_\nu(x)=\phi_\mu(Qx).
\]
Their zero-gradient regions are $C=[-c,c]^d$ and $Q^{\mathsf T}C$.
Both lie in the source cube, because
$(1-\varepsilon)(\cos\theta+\sin\theta)<1$, and both have mass
$p_\varepsilon=(1-\varepsilon)^d$.

Signed coordinate permutations show that the nonzero gradient atoms
of $\phi_\mu$ have equal mass. For $\phi_\nu$, the same conclusion follows from
quarter-turns in each coordinate plane and permutations of the
planes: these transformations preserve the source and commute with
$Q$, and act transitively on the nonzero gradient atoms. The target
measures $\mu=(\nabla\phi_\mu)_\#\rho$ and
$\nu=(\nabla\phi_\nu)_\#\rho$ are therefore
\begin{align*}
 \mu&=p_\varepsilon\delta_0+
   \frac{1-p_\varepsilon}{2d}\sum_{j=1}^d(\delta_{e_j}+\delta_{-e_j}),\\
 \nu&=p_\varepsilon\delta_0+
   \frac{1-p_\varepsilon}{2d}\sum_{j=1}^d
          (\delta_{Q^{\mathsf T}e_j}+\delta_{-Q^{\mathsf T}e_j}).
\end{align*}
Each atom of $\mu$ is paired with its nearest atom of $\nu$, and
these pairings respect the masses. The resulting coupling is optimal, so
\[
 W_2(\mu,\nu)^2=2(1-p_\varepsilon)(1-\cos\theta).
\]
On $C\mathbin{\triangle}Q^{\mathsf T}C$, one gradient is zero and the
other has norm one. Thus
\[
 \norm{T_\mu-T_\nu}_{L^2(\rho)}^2
 \ge2p_\varepsilon\bigl(1-r(\theta)^{d/2}\bigr),\qquad
 r(\theta)=\frac{2}{1+\sin\theta+\cos\theta}.
\]
Here $r(\theta)$ is the relative area of intersection of a square
and its rotation. To compute it, the removed part consists of four
corner triangles, each with legs
$c(\sin\theta+\cos\theta-1)/\sin\theta$ and
$c(\sin\theta+\cos\theta-1)/\cos\theta$. Dividing their total area
by $4c^2$ gives
\[
 r(\theta)=1-
 \frac{(\sin\theta+\cos\theta-1)^2}{2\sin\theta\cos\theta}.
\]
The $d$-dimensional overlap is the product of $d/2$ planar overlaps.
For each fixed even $d$, as $\varepsilon\downarrow0$,
\[
 W_2(\mu,\nu)^2\sim\frac d4\varepsilon^3,
 \qquad
 2p_\varepsilon\bigl(1-r(\theta)^{d/2}\bigr)
 \sim\frac d4\varepsilon.
\]
Taking $\varepsilon=(n+2)^{-1}$ defines targets $\mu_n$, $\nu_n$ with
\[
 W_2(\mu_n,\nu_n)^2\sim\frac{d}{4n^3},\qquad
 \norm{T_{\mu_n}-T_{\nu_n}}_{L^2(\rho)}^2
 \ge(1-o(1))\frac{d}{4n}.
\]
Since $k_d=d$ for even $d$, this proves the claimed estimates in this case.
For odd $d\ge3$, extend the construction in dimension $d-1$ by
making the potentials independent of the last coordinate. The maps
and target distances retain their values, so the claimed estimates
hold with $k_d=d-1$.
\end{proof}

\appendix
\section{Isotropic BV bound}
\label{sec:explicit-bv}

\begin{proof}[Proof of Proposition~\ref{prop:isotropic-bv}]
Use the normalization~\eqref{eq:source} and, as in
Lemma~\ref{lem:weighted-bv}, write
\[
 K_t=\{V\le t\},\qquad q(t)=\lvert K_t\rvert,\qquad
 p(t)=\Per(K_t).
\]
Each $K_t$, $t>0$, is a convex body. We first show that
\[
 \overline B(0,1/3)\subset K_{3d}.
\]
By the Brunn--Minkowski inequality, $h=q^{1/d}$ is nonnegative,
nondecreasing and concave. Hence $h(t)/t$ and $h'(t)$ are
nonincreasing, and so is
\[
 \frac{q'(t)}{t^{d-1}}
 =d\left(\frac{h(t)}t\right)^{d-1}h'(t),
\]
with derivatives taken almost everywhere on $(0,\infty)$.
For $X\sim\rho$, the law of $V(X)$ has an atom $q(0)/Z$ at zero
and density $Z^{-1}e^{-t}q'(t)$ on $(0,\infty)$. The ratio of this
density to the gamma density $t^{d-1}e^{-t}/(d-1)!$ is therefore
nonincreasing.
Comparing the two probability laws gives, for every $s>0$,
\[
 \rho\{V>s\}\le\frac1{(d-1)!}\int_s^\infty t^{d-1}e^{-t}\dd t.
\]
At $s=3d$, the right-hand side is $7e^{-6}$ for $d=2$ and
$(101/2)e^{-9}$ for $d=3$. For $d\ge4$, the exponential Markov
inequality bounds it by $(3e^{-2})^d\le81e^{-8}$. Each of these
bounds is less than $e^{-1}-1/3$, so
\[
 \rho(K_{3d}^c)<e^{-1}-1/3.
\]
For every unit vector $n$, the random variable $Y=n\cdot X$ is
isotropic and log-concave on $\R$. Its density is bounded by $1$,
and Gr\"unbaum's inequality gives $\Prob(Y>0)\ge e^{-1}$;
see~\cite[Lemmas~2.4 and~2.5]{MarsigliettiPandeyStatisticalDistances},
letting the concavity parameter tend to zero there. Consequently,
\[
 \Prob(Y>1/3)\ge e^{-1}-1/3.
\]
If $K_{3d}$ failed to contain $\overline B(0,1/3)$, separation
would give a unit vector $n$ with
$K_{3d}\subset\{x:n\cdot x\le1/3\}$, contradicting the two
preceding bounds. This proves the inclusion.

For $t>0$ and $\varepsilon>0$, convexity now gives
\[
 K_t+\varepsilon\overline B(0,1/3)
 \subset(1+\varepsilon)
 K_{(t+3d\varepsilon)/(1+\varepsilon)}.
\]
Comparing the right derivatives of the volumes at
$\varepsilon=0$ yields
\[
 \frac13p(t)\le d\,q(t)+(3d-t)q'(t)
 \qquad\text{for almost every }t>0.
\]
By log-concavity, $q$ grows at most polynomially.
The coarea formula from Lemma~\ref{lem:weighted-bv} and integration
by parts therefore give
\[
 \begin{aligned}
 \frac13\lvert Df\rvert(\R^d)
 &\le\frac1Z\int_0^\infty e^{-t}
       \bigl(d\,q(t)+(3d-t)q'(t)\bigr)\dd t\\
 &=4d-\int V\dd\rho-\frac{3d\,q(0)}Z
 \le4d.
 \end{aligned}
\]
This proves the proposition.
\end{proof}

\printbibliography

\end{document}